\documentclass[12pt,a4paper]{article}
\usepackage{amssymb}

\usepackage{graphicx}
\usepackage{amsmath}
\usepackage{amsmath,amssymb}%

\usepackage{amsfonts}
\usepackage{amsthm,amscd}
\usepackage[english]{babel}
\usepackage{amsmath}
\usepackage{amssymb}
\usepackage{amstext}

\newtheorem{theorem}{Theorem}

\newtheorem{example}[theorem]{Example}

\begin{document}

\title{Ergodic Theory in Classical and Bayesian inference}

\author{Artur Oscar Lopes \\ Inst. Mat. Est. -  UFRGS \\ Av. Bento Goncalves 9500 - Porto Alegre  - 91500-000 -  Brazil}

\date{July 16, 2024}

\date{\today}

\maketitle

\begin{abstract}We begin by presenting the mathematical rationale underlying classical deductive inference. We then introduce the foundational ideas of the Bayesian inference framework. Results lying at the interface of Statistics and Ergodic Theory are outlined, providing a theoretical framework applicable to the prediction and analysis of real-world phenomena from random data. This text is expository in nature — no new results are presented; rather, recently published results are described in a didactic manner. Throughout, we work with Hölder equilibrium measures, which encompass a substantially more general class of processes than i.i.d. ones.

\end{abstract}


\vspace{2mm}


\smallskip

email: arturoscar.lopes@gmail.com

  \section{Inductive Inference}
  
  This short note is intended to explain to mathematicians some applications of Ergodic Theory
to some specific problems in  Classical and  Bayesian inference.  We believe there is a need for an exposition that deals with inference in a didactic way and without major technicalities.    In principle, the Bayesian framework is opposed to the frequentist point of view. More explanations are needed on this issue.

In Section \ref{Clas} we will provide some simple examples for the reader that is not familiar with the topic of inference.

 At the end of the paper, we will consider inference results for H\"older equilibrium probabilities (see \cite{PP}), which is a class much more general than the classical i.i.d. processes; any shift invariant probability can be weakly approximated by a H\"older equilibrium probability (see \cite{L3} or \cite{GKLM}).

Recently the authors
 A. Nobel, K. McGoff, S.Mukherjee, and N.S. Pillai presented several interesting results in the interplay of Statistics and Ergodic Theory.
We will describe in a didactic and synthetic way another line of investigation: we highlight and put in context  some new results by M. Denker, H. H. Ferreira, M. Kesseb\"ohmer, A. O. Lopes, S. R. C. Lopes, J. Mengue, and  P. Varandas

In the frequentist setting, one assumes that the probability of an event is proportional to the times that the event occurs in a sequence of observations of random data; in this way, the Ergodic point of view is quite pertinent. The section \ref{Clas} will describe  basic results for classical inference under such a point of view.

Frequentists say the data are random and the parameters are fixed; on the other hand,  Bayesians say   the data are fixed, and the parameters are random.

In real-world problems, sometimes a satisfactory  sequence of data is available for a frequentist analysis, and in some cases, not. The choice of a purely Bayesian, or purely frequentist analysis depends very  much on the type of the problem, and the conviction of the investigator; it also depends on the relevance or not of the available data (a subjective matter) for the formulation of the Statistical model in question.

In the frequentist framework, typically, there is a fixed probability that generates a time series of random data. On the other hand, in the Bayesian setting, one can consider randomness in the probabilities that are relevant to the problem under consideration, or that the time series that eventually exists is small, or nonexistent,  and thus of little relevance.
In several cases, some subjectivity in the formulation of the model is necessary.
\medskip

From section 6.1 in the  book \cite{AR} we get the claim:

``Although the two approaches are conceptually different, they are nevertheless
not {\it complete strangers} to each other and may benefit from cooperation. Frequentist analysis of Bayesian procedures can be very useful for better understanding and
validating their results.''

\smallskip

In the present work, we consider some aspects of the Bayesian framework where one can take advantage of Ergodic Theory.
The final results will be  produced by mixing these two apparently opposed points of view. This will be considered in Section \ref{Bay}.

\smallskip

What is Inference in general terms?

 Quoting Chapter 1 in \cite{Cat}

 {\it The process of drawing conclusions from available information is called inference. When the available information is sufficient to make unequivocal, unique
assessments of truth, we speak of making deductions: based on a certain
piece of information, we deduce that a certain proposition is true. The method
of reasoning leading to deductive inferences is called logic. Situations where the
available information is insufficient to reach such certainty lie outside the realm
of logic. In these cases, we speak of doing inductive inference, and the methods
deployed are those of probability theory and entropic inference.}

\medskip

Given a certain problem to be analyzed, the use of the frequentist or Bayesian point of view may be a matter of dispute among different researchers. The use of each of the two options and their consequent predictions in comparison with the real world may show in each case whether it is more relevant to use one or the other.

Note that in the formulation of a statistical model, even when there are preliminary historical data, their quantity can be understood as sufficient (a subjective question) for the corresponding analysis or not. In the latter case, a Bayesian researcher often uses simulations with random data, obtained via a computer, to calibrate, or validate,  the model, and/or, to prove its efficiency (see \cite{Ka}). In this case, it can be said that there is a certain frequentist remnant. In some cases, such as when there is a long dependency in the data, obtaining long time series via random simulations on a  computer may be time-consuming (see, for instance, \cite{SL}, \cite{Olb},  \cite{Samo}, and \cite{ReLo}).

We do not intend to exhaust all the possible topics covered by Bayesian statistics, but only to describe in an intelligible mathematical way some aspects that may allow the reader to understand some of its essential ingredients (see beginning of Section \ref{Bay}).  For most of the models in Bayesian Statistics, the Ergodic point of view is not so relevant (see for instance, \cite{Wa}, \cite{AR}, and \cite{Bor}).

Some results relating Ergodic Theory and  Statistics are    \cite{Cha1}, \cite{D1}, \cite{DK1}, \cite{DK2}, \cite{DK5}, \cite{Denk1}, \cite{Denk2}, \cite{FLL},  \cite{Chan}, \cite{Leu}, \cite{Ent}, \cite{LLV}, \cite{LM7},  \cite{McG}, \cite{Nobel1}, \cite{Ka1}, \cite{LL1}, \cite{LL2} and \cite{Coll}.

In \cite{Cat} the author explores the point of view of deriving the main postulates of Physics from inference.

\medskip

\section{Classical Inference} \label{Clas}

An appropriate form of describing states of  partial belief (or uncertainty) is via  probabilities $q$ on a certain set $X$. A quite natural strategy is the following: the
 state of  partial belief $q^0$ should be updated to a new state of  partial belief $p^0$ when new information is available.
 The mathematical procedure to get the posterior probability $p^0$ is based on a mathematical formalism that takes into account the prior probability $q^0$ and also
 a function $l$, which is usually called the loss function (sometimes called risk function). The loss function may change according to the specific problem one is considering.
 We say that $p^0$ was obtained via an inference process obtained  from $q^0$ and $l$.

 The new information may arrive in the form of a random time series, but not necessarily in this form. In some cases, this information has to be processed via a minimization  procedure
 involving   probabilities.

 The states of partial belief may be indexed by a set $\Theta$ of possible probabilities.
 This way, we get a family of probabilities $p(\theta)$ for $\theta \in \Theta$. The probabilities $q^0$ and $p^0$ are among this set of possible probabilities.

 The randomness of the probability $q^0$ may be caused by lack of information, and inductive inference is  commonly used  in Information theory.

A simple example of inference can be described in the following way. Consider the case where $X=\{1,2,..,d\}$ and
$$\Theta= \{p=(p_1,p_2,..,p_d)\,|\, \sum_{j=1}^d  p_j  =1, p_j \geq 0, j=1,2,....,d\}.$$

\medskip

 Probabilities on $\Theta$ are denote by $p(\theta).$  The parameter $\theta$ is indexed by $(p_1,p_2,..,p_d).$

 We denote by $\Omega$ the symbolic space $X^\mathbb{N}$.
 Each probability $p(\theta)$ is associated with i.i.d. Stochastic Process, which means a probability $\mu_{\theta}$ on $\Omega.$

We consider a {\bf special point $\theta_0 \in \Theta$ which will be responsible for producing random sequences} $y = (y_1,y_2,...,y_n,...)\in \Omega$. The associated probability on $\Omega$ is denoted  by $\mu_{\theta_0}$, which in principle we don't know.
 However, we can observe finite random samples over time.
 Given a sample  sequence $y= (y_1,y_2,...,y_n,...)\in \Omega $,  we denote $y^n$ the finite word $y^n= (y_1,y_2,...,y_n) \in X^n$. Moreover, we denote by $C_n(y)$ the cylinder $\overline{y_1,y_2,...,y_n}=\overline{y^n}$.

A very useful transformation in the study of probability measures in the symbolic space $\Omega$ is the shift transformation $\sigma:\Omega \to \Omega$. We set
$$\sigma(y_1,y_2,...y_n,...)= (y_2,y_3,...,y_n,...).$$

An Stochastic Process with values on $X$ determine a probability $P$ on $\Omega$. A stochastic process being stationary corresponds to the property that for any cylinder $C_n(y)$, it holds that
$P( \sigma^{-1} (C_n(y)  ))=P ( C_n(y) ).$ If this property is true we say that $P$ is invariant for the action of the shift.
\medskip

 In most cases, for applications  in the  real world, a fixed  random source  is the object of interest in the problem under consideration, which is not precisely  known but generates random samples.

The finite sample  $y^n= (y_1,y_2,...,y_n)$ corresponds to the new information at time $n$.

 We don't know which $\theta_0$ produced the sample $y$, but if we have to bet, it would be better to choose the $\hat{\theta}= \hat{\theta}(y)$ such that
$$  \mu_{\hat{\theta}}( C_n(y) ) \geq  \mu_{\theta} ( C_n(y) ),$$
when compared to other possible $\theta$. In this case, it is natural to say that
$\mu_{\hat{\theta}}$ is more
suited to the data   $y^n$ than $\mu_{\theta}$ (see Section 1.2 in \cite{AR}).  We will elaborate on this claim.

The likelihood ratio
$$ \frac{ \mu_{\hat{\theta}}( C_n(y) )}{ \mu_{\theta_0} ( C_n(y) )}$$
shows the strength of the evidence in favor of $\hat{\theta}= \hat{\theta}(y) $ against $\theta_0$.

In this way is natural to try to maximize
$$ \log \mu_{\hat{\theta}}( C_n(y) ) - \log \mu_{\theta_0} ( C_n(y) ).$$

Then, it is natural to  consider  a minimization problem associated with the {\em loss functions}
$\ell_n : \Theta \times \Omega  \to \mathbb R$, $n\in \mathbb{N}$, given by
\begin{equation} \label{Chalk}
\ell_n(\theta,y)=- \log \Big(\frac{ \mu_{\theta}( C_n(y) ) }{ \mu_{\theta_0} ( C_n(y) )} \,\, \Big)=  \log \Big(\frac{ \mu_{\theta_0}( C_n(y) ) }{ \mu_{\theta} ( C_n(y) )} \,\, \Big)
\end{equation}

This loss function is of  the log-likelihood form. We refer the reader to Chapter 2, Section 16 in \cite{Bor} for general results on the maximum likelihood method.

At first, one can ask if there exists a conceptual issue with the above reasoning because we  claim we do not know $\theta_0$, but $l_n$ encapsulates the information of $\theta_0$. We will show that this will not be a problem  (see \eqref{Cha4}).

 For a given $n \in \mathbb{N}$,  and $y$ chosen at random with respect to $\theta_0$, we are interested in the value $\theta^n=\theta^n_y\in \Theta$ minimizing $\theta \to \ell_n(\theta,y)$. In this case $p(\theta^n)$ was obtained by inductive inference from $l$ and $p(\theta_0).$

 We would like to show that for $y$ chosen at random with respect to $\theta_0$,  we get that $\theta^n \to \theta_0$ when $n \to \infty$.
 The bottom line is: from the loss function $l$, by getting more and more information, we will be able to determine the source responsible for randomness.
 \smallskip

\begin{example}
Set $\theta_0=\tilde{p}= (1/3,1/3,1/3)$, and for fixed  $n>0$, consider a  random sample $y=(y_1,y_2,...,y_n)$ of size $n$, obtained from the independent probability $P$ on $\{1,2,3\}^\mathbb{N}$ associated to $\tilde{p}$. 
The observer knows in advance that the source is produced by an independent probability $P$
but does not know the exact $\theta=(q_1,q_2,q_3)$, $q_1+q_2+q_3=1$, that produces it. From the information of the sample $y$ he wants to obtain the best approximate values $\hat{\theta}=(\hat{q}_1,\hat{q}_2,\hat{q}_3).$ The probability $\tilde{p}$ is unknown to the observer. In this example we follow the frequentist point of view.
 
As we obtain a sample of larger size $n$, one can obtain a better approximation of the true values of the unknown $\theta=(q_1,q_2,q_3)$ that determine the probability $P$. Given a sample $y$, what is the best guess? We believe this simple example will illustrate the effectiveness of the maximum likelihood estimator procedure for obtaining a good guess $\hat{\theta}=\hat{\theta}_n=(\hat{q}_1,\hat{q}_2,\hat{q}_3).$
We denote by $C_n(y)$ the cylinder set $\overline{y_1,y_2,...,y_n}\in \{1,2,3\}^\mathbb{N}$. The function $g_n$ will be defined as the maximum likelihood estimator
$$g_n(p)=  \frac{ 1 }{ n}  \log \left(\frac{p (   C_n (y))}{\tilde{p} (   C_n (y) )} \right) .$$
 For fixed $n$, one is interested in the probability $p_ n=(p^1_n,p^2_n,p^3_n)$ maximizing $g_n$. The loss function is $-\,g_n$.
Assume the string $y$ has $n_j$ occurrences of each symbol $j=1,2,3.$ In this case $n_1+n_2+n_3=n$. Using Lagrange multipliers one can show 
that the probability $p_n$ maximizing the above expression is such that
\begin{equation} \label{toaq2} p^j_n= \frac{n_j}{n}, \quad j=1,2,3.
\end{equation}
We elaborate on this: there are two constraints $p_1 +p_2 + p_3=1$ and $n_1 +n_2 + n_3=n$. The gradient of $G(p_1,p_2,p_3)=  p_1 +p_2 + p_3$ is $(1,1,1).$
The gradient of 
\begin{equation} \label{toaq5}F(p_1,p_2,p_3) = \log (p_1)^{n_1}  \,(p_2)^{n_2} \,(p_3)^{n_3}- \log 3^{-n}
\end{equation}
is
$$\nabla F (p_1,p_2,p_3) =\left(  \frac{n_1}{p_1}, \frac{n_2}{p_2}, \frac{n_3}{p_3}\right),\,\, \text{and we assume it is equal to}\,\, (\lambda, \lambda,\lambda).$$
That is, the Lagrange multipliers method produce
\begin{equation} \label{toaq1}\lambda=  \frac{n_1}{p_1}, \quad \lambda=\frac{n_2}{p_2}, \quad \lambda=\frac{n_3}{p_3}.\end{equation}
From the above we derive 
\begin{equation} \label{toaq3}\lambda=\lambda (p_1+p_2+p_3)=n_1+n_2+n_3=n.
\end{equation}
Therefore, from \eqref{toaq1} and \eqref{toaq3} it follows that \eqref{toaq2} holds.
\medskip

One can show that $p_n\to \tilde{p}$ as $n\to \infty$. Indeed, for fixed $j$,
denote $I_{\overline{j}}: \{1,2,3\}^\mathbb{N} \to \mathbb{R}$ the indicator function of the cylinder $\overline{j}\subset  \{1,2,3\}^\mathbb{N}$.
Note that, by the Law of Large Numbers (or by the Birkhoff Ergodic Theorem),
$$\frac{n_j}{n}= \frac{\sum_{k=0}^{n-1}\,I_{\overline{j}} (\sigma^k (y))\,}{ n}\to \int I_{\overline{j}}\,d P= P( \overline{j} )=1/3.$$ 
The sample is finite: suppose $y=(2,1,2,1,3,2,1,3,2,3)$, where $n= 10$. The guess is
$$ \hat{\theta}=(\hat{q}_1,\hat{q}_2,\hat{q}_3)=(3/10, 4/10,3/10). $$
The reasoning above would work for a further problem with other values $\theta_0=\tilde{p}\neq (1/3,1/3,1/3)$, producing different samples (here the ergodic property play its role), and deriving others $ \hat{\theta}$ .

\,\,\,\,\,\,\,\,\,\,\,\,\,\,\,\,\,\,\,\,\,\,\,\,\,\,\,\,\,\,\,\,\,\,\,\,\,\,\,\,\,\,\,\,\,\,\,\,\,\,\,\,\,\,\,\,\,\,\,\,\,\,\,\,\,\,\,\,\,\,\,\,\,\,\,\,\,\,\,\,\,\,\,\,\,\,\,\,\,\,\,\,\,\,\,\,\,\,\,\,\,\,\,\,\,\,\,\,\,\,\,\,\,\,\,\,\,\,\,\,\,\,\,\,\,\,\,\,\,\,\,\,\,\,\,\,\,\,\,\,\,\,\,\,\,\,\,\,\,\,\,\,\,\,\,\,\,\,\,\,\,\,\,\,\,\,\,\,\,\,\,\,\,\,\,\,\,\,\,\,\,\,$\diamond$

\end{example}

 Following \cite{Cha1}, given the ergodic probability measures  $\mu_{\theta_0}$   and $\mu_\theta$, the relative entropy $ h(\mu_{\theta} \mid \mu_{\theta_0})$ is given by the limit
$$ h(\mu_{\theta_0} \mid \mu_{\theta}):=\lim_{n \to \infty} \frac{ 1 }{ n}  \log \Big(\frac{ \mu_{\theta_0} (   C_n (y))}{\mu_{\theta}  (   C_n (y) )} \Big)
=$$
$$ \sum_{r=1}^d p(\theta_0)_j \log p(\theta_0)_j - \sum_{r=1}^d p(\theta)_j \log p(\theta_0)_j =$$
\begin{align} \label{Chacha5}
  \int   p(\theta_0) \log  p(\theta_0) -\int   p(\theta) \log  p(\theta_0)\geq 0,
\end{align}
which exists and is non-negative, for $\mu_{\theta_0}$-almost every $y=(y_1,y_2,...,y_n,..)\in \Omega$. The minimum value of \eqref{Chacha5} is zero, and it is  attained just when $\theta=\theta_0.$

Note that $h(\mu_{\theta_0} \mid \mu_{\theta})=0$, if and only if, $\theta=\theta_0$. In the present case, if $h(\mu_{\theta_0} \mid \mu_{\theta_n})\to 0$, then $\theta_n \to 0$, when $n \to \infty$.

Given $\theta_1,\theta_2$, if
$$ \int   p(\theta_0) \log  p(\theta_0) -\int   p(\theta_1) \log  p(\theta_0)>$$
\begin{align} \label{Cha1}
    \int   p(\theta_0) \log  p(\theta_0) -\int   p(\theta_2) \log  p(\theta_0),
\end{align}
then $\mu_{\theta_2}$ is preferred, when compared with  $\mu_{\theta_1}$.

Note that \eqref{Cha1} is equivalent to
\begin{align} \label{Cha2}
  -\int   p(\theta_1) \log  p(\theta_0)>    -\int   p(\theta_2) \log  p(\theta_0),
\end{align}
and then, given the finite  random sample $y^n=(y_1,y_2,...,y_n)$, it is natural to say, taking into account the limit, that  for large $n$,  $\mu_{\theta_2}$ is preferred when compared with  $\mu_{\theta_1}$,
if
$$ \mu_{\theta_2}  (   C_n (y) ) \geq \mu_{\theta_1}  (   C_n (y) ).$$
Then, given the random data $y^n=(y_1,y_2,...,y_n)$, from the inference point of view, we should   look for the $\theta$ maximizing
\begin{equation} \label{Cha4} \mu_{\theta}  (   C_n (y) )=  \mu_{\theta} (\overline{y_1,y_2,...,y_n} ).
\end{equation}
A natural question is: in the long run, taking large $n$ in the above  procedure, in the limit, will  we localize $\theta_0$?
We will see that the answer is yes.

In the maximum-likelihood estimator procedure (see definition 1 in Chapter II.16 on page 69 in \cite{Bor}), for fixed $\theta_0$,  we want  to maximize in $\theta$ the value
\begin{align} \label{Chachal}
   \int   p(\theta) \log  p(\theta_0) -  \int   p(\theta_0) \log  p(\theta_0)\leq 0,
\end{align}
which is equivalent to minimize   in $\theta$ the expression \eqref{Chacha5}.

The main question is: in the long run, inggetting larger and larger samples, can we identify $\mu_{\theta_0}$?

 Suppose  $E\subset \Theta$ is a closed set such that $\theta_0\notin E$, then, there exists an $\alpha>0$, such that,
\begin{equation} \label{fert}\inf \{  \sum_{r=1}^d p(\theta)_j \log p(\theta_0)_j - \sum_{r=1}^d p(\theta_0)_j \log p(\theta_0)_j \,|\, \theta \in E\,\}  <- \alpha <0.
\end{equation}

 %\eqref{fert} is equivalent to say that
% $$E=\{\theta \,|\,   \int   p(\theta) \log  p(\theta_0) -\int   p(\theta_0) \log  p(\theta_0)\leq - \alpha\}.$$

 We can choose as an example  $E$ of the form $\{\theta \,| \, |\theta - \theta_0 |>\epsilon\}$, for small $\epsilon$.

  Then, for $\theta \in E$, and for $\mu_{\theta_0}$-almost every $y=(y_1,y_2,...,y_n,..)\in \Omega$, we get  from \eqref{Chacha5}, that for large $n$
 \begin{equation} \label{fert1}- \log (\frac{ \mu_{\theta} (   C_n (y))}{\mu_{\theta_0}  (   C_n (y) )} ) \geq   \alpha\, n.
 \end{equation}

  Suppose  $F$ is a closed set of the form
 $$\{\theta \,|\,   - \frac{\alpha}{2}\leq  \int   p(\theta) \log  p(\theta_0) -\int   p(\theta_0) \log  p(\theta_0)\leq - \frac{\alpha}{4}\}$$ not containing $\theta_0$. Then, for $\theta \in F$, and for $\mu_{\theta_0}$-almost every $y=(y_1,y_2,...,y_n,..)\in \Omega$, we get for large $n$
 \begin{equation} \label{fert4}\frac{\alpha}{2}\, n  \geq - \log (\frac{ \mu_{\theta} (   C_n (y))}{\mu_{\theta_0}  (   C_n (y) )})  \geq  \frac{\alpha}{4}\, n.
 \end{equation}

 Note that $F \cap  E = \emptyset.$
  Therefore, given $n$, for  $\mu_{\theta_0}$-almost every $y=(y_1,y_2,...,y_n,..)\in \Omega$, when looking for a $\theta$  minimizing
 $$ \ell_n(\theta,y)=- \log \Big(\frac{ \mu_{\theta} \,( C_n(y) }{ \mu_{\theta_0}\,( C_n(y) ) } \, \Big),$$
 the optimal $\theta$ will not be in the set $E$ for sure. Indeed, taking large samples $y^n$, in the long range, a $\theta_2$ in $F$ is preferred when compared with a $\theta_1$ in $E$ (see \eqref{Cha1}); this is so because, $\alpha\, n > \frac{\alpha}{2} \, n.$

 In this way, the minimizing solution $\theta^n=\theta^n_y$, for  the function $\theta \to \ell_n(\theta,y)$, for large $n$, will not be in the set $E$.

 A full proof of the results described above follows from the reasoning of \cite{LLV}.

The  relative entropy is also known  as  Kullback-Leibler divergence (see \cite{FLL} and \cite{LM7}).

 \section{Bayesian Inference} \label{Bay}

 We present here a particular point of  view concerning a certain class of problems in the topic, which is very broad; other frameworks are also relevant and interesting, but not covered here.

 In the same way as in the last section,
the states of partial belief are indexed by a set $\Theta=[0,1]^n$.  Now the elements in $\Theta$  are probabilities  $\mu$ on $\Omega=\{1,2,...,d\}^\mathbb{N}$.
 In this way, we get a family of probabilities $\mu_\theta,$ for $\theta \in \Theta$.  But now, in the Bayesian framework, we assume that the parameters $\theta$ are in
 a state of partial belief which is described by a probability $Q_0$ in $\Theta.$
This probability $Q_0$ may be associated with complete a priori ignorance, or alternatively, describe some subjective belief of those who wish to apply the method.

I.  In one possible setting, there exists a special unknown   parameter $\theta_0$, such that $\mu_{\theta_0}$ produces the relevant samples, in a similar way as in the last section;  but, now we set an initial distribution of probability given  by a certain $Q_0$ on $\Theta$, describing the uncertainty on $\theta_0$.

II.  In another possible setting, each $\mu_\theta, \theta \in \Theta$, produces samples. In this case, the role of the probability $Q_0$ is to average the randomness of the samples obtained from the random family $\mu_\theta, \theta \in \Theta$.

  In any case, in principle, $Q_0$ describes fuzzy information on probabilities. From different  types of inference procedures, given a loss function $l$, not necessarily of log-likelihood type, one can get a posterior probability  $Q_1$ in $\Theta$, from the prior probability  $Q_0.$

 The loss function $l$ can have different forms depending on the problem of  interest, and $l$ may incorporate, or not, some subjective belief of those who wish to apply the method. This applies to both problems: types I and II.

 One should consider two more  alternatives that apply to cases I and II: a) whoever wishes to apply the method has a pertinent finite random sample of arbitrarily long size at their disposal; b) the  item a) is not fulfilled.

 In the first case a),   the frequentist point of view could be incorporated into the Bayesian framework. Then, taking into account the sequence of  finite random data $y^n = (y_1,y_2,...,y_n) $, $n \in \mathbb{N},$ and a family of loss functions $l_n$, we can obtain a family of probabilities $Q_n$ in $\Theta$. In this case, the information of $Q_n$, $n \in \mathbb{N},$ is less uncertain than the one given by $Q_0$. We are interested in the eventual limit of such a family $Q_n$ on $\Theta$, when $n \to \infty$.

 In case b), using an specific $l$, one gets the posterior $Q_1$ from the prior $Q_0$, and this should satisfy (by particular reasons inherent to the problem)  the properties aimed at by the one who applies the method. In this case, in the establishment of the Statistical model, a greater degree of subjectivity is inherent in its formulation. The choice of $l$ needs to incorporate the appropriate relevance that is natural to the problem under consideration

\smallskip

Here we choose to describe the case where we assume I) and a).  Subjectivity in the formulation of this model is not  necessary. Other cases   are of interest, of course, and we refer the reader to \cite{AR} and \cite{Bor}.

\smallskip

We point out that the Bayes formalism, and the Bayes rule, can be deduced via an  Optimal Information procedure, as described in \cite{Zel};  the topic was the subject of a more in-depth study, and considered in a more general setting in \cite{Ent}.

Assume that  the probability $\mu_{\theta}$ on $\Omega= \{1,2,...,d\}^\mathbb{N},$ is  the equilibrium probability for a H\"older potential  $A_\theta :\Omega \to \mathbb{R}$ (see \cite{PP}) indexed by a parameter $\theta \in \Theta=[0,1]^k$  (see \cite{Nobel1}, \cite{Denk1} and \cite{LLV}). We will mention some new inference results concerning such a  family.

 For a special fixed parameter $\theta_0\in \Theta$ we will get  data from the probability $\mu_{\theta_0}$, which is the equilibrium probability for the fixed H\"older potential $A_{\theta_0}$

For convenience, one can assume that the support  of $Q_0$ is the whole space $\Theta$.
As an example, $Q_0$ could be the uniform probability on the simplex $\Theta$ (but we do not have to assume that).
This could be seen as a form of total ignorance. The initial choice of the probability $Q_0$ will not be relevant in our setting.
Can we identify from the random data, via an iteration procedure,  the probability $\mu_{\theta_0}$ among the others in $\Theta$?

We will consider once more the log-likelihood point of view.

We  assume that the {\em loss functions}
$\ell_n : \Theta \times \Omega \to \mathbb R$, $n\in \mathbb{N}$, are given by
%

%\marginpar{\tiny this loss function depends on $\theta_0$ BUT independs on $\theta$; is this desired?}
%
\begin{equation} \label{Chacha}
\ell_n (\theta, y) = - \log \Big(\frac{ \mu_{\theta}( C_n(y) ) }{ \mu_{\theta_0} ( C_n(y) )} \,\, \Big)\,
 \end{equation}

We set
\begin{align*} 
Z_n(y) & =  \int_\Theta  \text{exp}^{ - \, \ell_n (\theta,y)} \,\, d Q_0 (\theta)
   \\
 &
  = \int_\Theta \,  \frac{\mu_{\theta}  \,(\,C_n(y) \,)}{\mu_{\theta_0} \,(\,C_n(y)\,) }     \,d Q_0 (\theta)
\end{align*}
for each $y\in Y$.

From \cite{LLV} we get
\begin{align}
\limsup_{n \to \infty} \frac{1}{n} \log Z_n(y) \,\, &  =  \limsup_{n \to \infty} \frac{1}{n} \log \int_\Theta  \frac{\mu_{\theta}
 \,(\,C_n(y) \,)}{\mu_{\theta_0} \,(\,C_n(y)\,) }   d Q_0 (\theta) \nonumber \\
 & \geqslant     \limsup_{n \to \infty} \frac{1}{n} \int_\Theta \log  \frac{\mu_{\theta}  (C_n(y) )}{\mu_{\theta_0} (C_n(y)) }
 %d \mu_{\theta_0}  (y)
 \, d \Pi_0 (\theta) \nonumber  \\
 & =  \int_\Theta  h(\mu_{\theta} \mid \mu_{\theta_0}  ) \, d Q_0 (\theta) >0
\end{align}
for $\mu_{\theta_0}$-almost every $y$.

 Given $n \in \mathbb{N}$ and $y\in \Omega$, chosen at random  according to $\mu_{\theta_0}$, we denote by $Q_n=Q_n^y$ the $n$-posterior probability on $\Theta$, which is given by
 \begin{equation} \label{frt}
E \to {Q_n^y (E)\,=\, \frac{ \int_E \mu_\theta (C_n (y) ) \,d Q_0 (\theta) }{ \int_\Theta \mu_\theta (C_n (y) ) \, d Q_0 (\theta)},}
\end{equation}
%\marginpar{\tiny is this expression OK? It seems it is the second one, am I wrong?}
on Borel sets $E \subset \Theta$.

In \cite{Nobel1} the authors call \eqref{frt} the  $n$-Gibbs posterior distribution.

Given an initial general $Q_0$, the proof that
$$\lim_{n \to \infty} Q_n^y   = \delta_{\theta_0},$$
for $\mu_{\theta_0}$-a.e. $y\in \Omega$,
is a non trivial generalization of the case where just i.i.d. probabilities are considered. Large deviation properties are used in the proof. We will elaborate on the main claim.

Note that  for any $n,y$ and $E$ we get
\begin{equation} \label{frt1}
Q_n^y (E)\,=\, \frac{ \int_E \,  \frac{\mu_{\theta}  \,(\,C_n(y) \,)}{\mu_{\theta_0} \,(\,C_n(y)\,) }     \,d Q_0 (\theta)}{ \int_\Theta \,  \frac{\mu_{\theta}  \,(\,C_n(y) \,)}{\mu_{\theta_0} \,(\,C_n(y)\,) }     \,d Q_0 (\theta)}= \frac{ \int_E  \text{exp}^{ - \, \ell_n (\theta, y)} \,\, d Q_0 (\theta) }{ Z_n (y)}.
\end{equation}

When $E$ is closed and $\theta_0 \notin E $ one can show from \eqref{fert1}  that there exists $\alpha>0$, such that
 \begin{equation} \label{fert2} \int_E \frac{ \mu_{\theta} (   C_n (y))}{\mu_{\theta_0}  (   C_n (y) )}  d  Q_0 (\theta) < \text{exp}^{-\alpha\, n}\, Q_0 (E).
 \end{equation}

In \cite{LLV}the following result    was shown :

\begin{theorem} \label{thm:main}
Given a Borel set  $E \subset \Theta$,
the limit
\begin{equation} \label{yrt}
      \lim_{n \to \infty} \frac{1}{n} \log Q_n^y (E)
      =  \lim_{n \to \infty} \frac{1}{n} \log    \frac{ \int_E \mu_\theta (C_n (y) ) d Q_0 (\theta) }{ \int_\Theta \mu_\theta (C_n (y) ) d Q_0 (\theta)}  \end{equation}
exists for $\mu_{\theta_0}$-almost every $y$. If $E$ is a closed set not containing $\theta_0$, the above limit is negative.
This means $  Q_n^y (E)$ goes to zero exponentially fast.

\smallskip

Moreover,
\begin{equation*}
      \lim_{n \to \infty} Q_n^y   = \delta_{\theta_0}, \quad \text{for $\mu_{\theta_0}$-a.e. $y\in \Omega$}.
  \end{equation*}
\end{theorem}

%There are no Competing Interests of any kind
\smallskip

%A. O. Lopes is partially supported by CNPq scholarship
%\smallskip

%No data of any nature was used in the work.

\smallskip
%This research did not receive any specific grant from funding agencies in the public, commercial, or not-for-profit sectors.

{}

\end{document}